\documentstyle[12pt]{article}
\def\k{{K(x,y)}}
\def\f{{\varphi}}
\def\e{{\varepsilon}}
\def\t{{T^\nu}}
\def\m{{M_n({\bf C})}}
\date{14 February 1995}
\author{I.~A.~Ikramov, \  F.~Sharipov}
\title{On the discrete spectrum of non-analytic matrix-valued Friedrichs
model}

\begin{document}
\maketitle
\begin{abstract}
We have found the sufficient conditions for the spectrum of
matrix-valued Friedrichs model to be finite.
\end{abstract}

\section{Introduction}
Consider a self-adjoint operator $H$ in the Hilbert space
$L^2(\t,{\bf C}^n)$ given by the formula
\begin{equation}\label{1}
(Hf)(x)=U(x)f(x)+\int_{\t}\k f(y)\,dy
\end{equation}
where $\t$ is a $\nu$-dimensional torus, ${\bf C}^n$ is the
$n$-dimensional complex Euclidean space, $L^2(\t,{\bf C}^n)$ is the complex
Hilbert space of the square-integrable (with respect to the norm) functions
on $\t$ taking values in ${\bf C}^n$, $U(x)$ and $\k$ are continuous
functions on $\t$ and $\t\times\t$ respectively with values in the space
$\m$ of complex $n\times n$-matrices, satisfying the following conditions:
\begin{equation}\label{2}
U^*(x)=U(x),\quad K^*(x,y)=K(y,x),
\end{equation}
where $^*$ denotes the adjoint matrix.

Operator of the form (\ref{1}) for the first time was treated by Friedrichs%
{}~\cite{1},~\cite{2} in the case $\nu=1$, $n=1$, $U(x)=x$ as a simple model
of the perturbation theory for continuous spectrum. His results were
further developed by O.~A.~Ladyzhenskaya and L.~D.~Faddeev~\cite{3},
L.~D.~Faddeev~\cite{4}, S.~N.~Lakaev~\cite{5},~\cite{6} etc.

In the case when $U(x)$ and $\k$ are matrix-valued analytic functions the
spectral properties of the operator $H$ were studied in~\cite{5},%
{}~\cite{6},~\cite{7},~\cite{8}. In particular the following theorem was
proved:

{\bf Theorem 1.}~\cite{6},~\cite{7}. {\it Let $U(x)=(\delta_{ij}u_j(x)),\
\k=(K_{ij}(x,y))$ where $\delta_{ij}$ is the Kronecker symbol and suppose
that the matrix elements $u_j(x)$, $K_{ij}(x,y)\ (i,j=1,2,\ldots,n)$ are
real-valued analytic functions on $\t$ and on $\t\times\t$ respectively. If
one of the following two conditions is valid then the operator $H$ defined
by (\ref{1}) has a finite number of eigenvalues (counted with
multiplicities) outside the essential spectrum.
\begin{enumerate}
\item
$\nu=1$ and for every $j=1,2,\ldots,n$ all critical points of the functions
$u_j(x)$ are isolated;
\item
$\nu\geq 1$ and for every $j=1,2,\ldots,n$ all critical points of the
functions $u_j(x)$ are non-degenerate.
\end{enumerate}
}

In the present paper we show the finiteness of the discrete spectrum of the
operator $H$ for a more wide class of functions $U(x)$ and $\k$.

\section{Main theorem}
Denote by $\sigma_{ess}(H)$ the essential spectrum of the
operator $H$ and denote by $\Delta(x,z)$ the determinant of the matrix
$U(x)-zE$, where $E$ is the unit matrix and $z\in{\bf C}$.

{\bf Lemma 2.}\ {\it The essential spectrum of $H$ is of the form
$$ \sigma_{ess}(H) = \cup_{x\in\t}\{z\in{\bf C}:\ \Delta(x,z)=0\}=
\cup_{i=1}^k [m_i,M_i] $$
with mutually non-intersecting segments $[m_i,M_i]\ (i=1,2,\ldots,k)$.}

Lemma 2 follows from the Weyl's theorem about the essential spectrum and
from the minimax principle (cf.~\cite{13}).

Denote now by $\Gamma$ the set
$$ \Gamma=\{m_1,m_2,\ldots,m_k,M_1,M_2,\ldots,M_k\} $$

{\bf Lemma 3.}\ {\it For any $z\in\Gamma$ the value $A=0$ is globally
extremal value for the continuous real-valued function
$\f_z(x)=\Delta(x,z)$ on the torus $\t$. }

{\bf Definition 1.}\ Let $\f(x)$ be a continuous real-valued function on
the torus $\t$. Extremal point $x^0\in\t$ of the function $\f(x)$ is
called a point of {\it finite multiplicity} if there exist such numbers
$m>0$,
$c>0$ and such neighbourhood $V(x^0)$ of the point $x^0$ that for any $x\in
V(x^0)$ the following inequality holds
\begin{equation}\label{3}
\vert\f(x)-\f(x^0)\vert\geq c\vert x-x^0\vert^m,\quad\mbox{where}\
\vert x-x^0\vert^2=\sum_{i=1}^n(x_i-x^0_i)^2.
\end{equation}
Else the extremal point $x^0$ is called a point of {\it infinite
multiplicity.}
The exact lower bound of the set of numbers $m>0$ satisfying the condition
(\ref{3}) is called {\it multiplicity}\/ of the extremal point $x^0$ and is
denoted by $m(x^0)$. Multiplicity of an extremal value $A$ of the function
$\f(x)$ is the sum of multiplicities of all extremal points from
the inverse image $\f^{-1}(A)$ of  $A$.

Denote by $C^{\alpha+0}(\t\times\t,\m)$ the space of the matrix-valued
functions $\k$ on $\t\times\t$ such that for any multiindex $\beta$ with
$\vert\beta\vert\leq [\alpha]$ the derivative $K^{(\beta)}(x,y)$ satisfies
the H\"older condition with index $\{\alpha\}+0$ where $[\alpha]$ is the
entire part of $\alpha$ and $\{\alpha\}=\alpha-[\alpha]$.

{\bf Theorem 4.}\ {\it Let $0<\mu<\infty$. Suppose that for any $z\in\Gamma$
number $A=0$ is the extremal value of the function $\f_z(x)=\Delta(x,z)$ of
the multiplicity $\leq\mu$. Let the function $\k$ belong to the class
$C^{2\mu-\nu/2+0}(\t\times\t,\m)$. Then the operator $H$ has only
a finite number of eigenvalues (counted with
multiplicities) outside the essential spectrum.}

{\bf Proof}\ of this theorem consists of the following lemmas.

{\bf Lemma 5.}\ {\it Let $z_0\in\Gamma$. If the matrix-valued function
$\Delta^{-1}(x,z_0)\,\k$ is square-integrable with respect to the norm on
$\t\times\t$, then there exists a positive number $\e=\e(z_0)$ such that
the operator $H$ defined by (\ref{1}) has only
a finite number of eigenvalues
(counted with
multiplicities)
in the set $(z_0-\e,z_0+\e)\setminus\sigma_{ess}(H)$.}

{\bf Proof.}\ By the Fredholm theorem it is sufficient to show that
$z_0$ is not a limit point of the discrete spectrum of the operator $H$.
Suppose that it is not so, i.e. there exists a sequence $\{z_n\}$ of
eigenvalues $(z_n\notin\sigma_{ess}(H))$, converging to $z_0$ and let $f_n$
be a normed eigenfunction of the operator $H$ corresponding to the
eigenvalue $z_n$, i.e. a solution of the equation
\begin{equation}\label{4}
(U(x)-z_nE)f_n(x)+\int_{\t}\k\,f_n(y)\,dy =0.
\end{equation}
Consider a sequence of operators
$$(\hat K (z_n)f_n)(x)=\int_{\t}(U(x)-z_nE)^{-1}\k\,f(y)\,dy,\quad
n=1,2,\ldots $$
By supposition $\hat K (z_n)$ is a compact operator and
\begin{equation}\label{5}
\lim_{n\to\infty}\hat K (z_n)=\hat K (z_0)
\end{equation}
in the uniform operator topology, hence the operator $\hat K (z_0)$ is
also compact.
Put $F=\{f_n:\ n=1,2,\ldots\}$. As the set $\hat K (z_0)F$ is precompact and
$$ f_n(x)=-\int_{\t}(U(x)-z_nE)^{-1}\k\,f_n(y)\,dy =-(\hat K (z_n)f_n)(x),\
n=1,2,\ldots $$
so by (\ref{5}) the set $F$ is also precompact. It contradicts the
orthonormality of the sequence $\{f_n\}$, and the lemma is proved.

{\bf Lemma 6.}\ {\it Let $B$ be a bounded self-adjoint operator in a
Hilbert space ${\cal H}$. If the essential spectrum of $B$ consists of
union of finite number of segments and if outside the essential spectrum
$B$ has a finite number of eigenvalues, then for any finite-dimensional
operator $K$ in ${\cal H}$ the operator $B+K$
has finite number of
eigenvalues
(counted with
multiplicities)
outside the essential spectrum of $B$.}

The lemma 6 can be easily proved using the Weyl theorem and the Fredholm
determinant.

{\bf Proof}\ of the theorem 4. If the conditions of the theorem 4 are
satisfied then for any $z\in\Gamma$ the function $\k$ can be represented
in the form $$ \k=K_1(x,y)+K_2(x,y) $$ so that the following conditions
are valid:
$$
\Vert\Delta^{-1}(x,z)K_1(x,y)\Vert\in L^2(\t\times\t),\ K^*_1(x,y)=
K_1(y,x),\ K^*_2(x,y)=K_2(y,x)
$$
and the integral operator with the kernel $K_2(x,y)$ in the space $L^2(\t,
{\bf C}^n)$ is finite-dimensional.

Now using lemmas 2,3,5 and 6 begin proving the theorem 4.

\section{Applications}
{\bf 1.}\
Suppose that the matrix-valued function $U(x)$ is analytic.
Then for any $z\in{\bf R}$ the function $\f_z(x)=\Delta(x,z)$ is
real-analytic on $\t$. It follows from the Lojasevitch inequality
(cf.~\cite{11}) that the isolated extremal point of a real-analytic
function on the torus $\t$ is an extremal point of finite multiplicity (cf.
definition 1). Henceforth by the theorem 4 we obtain the following theorem
generalizing the theorem 1.

{\bf Theorem 7.}\ Let the matrix-valued function $U(x)$ be analytic on $\t$
and let for any $z\in\Gamma$ the set $\f^{-1}_z(0)$ is finite (where
$\f_z(x)=\Delta(x,z)$). Then there exists a positive number $s>0$ such
that for any matrix-valued function $\k$ from the class
$C^{s+0}(\t\times\t,\m)$ the operator $H$
has finite number of
eigenvalues
(counted with
multiplicities)
outside the essential spectrum.

The following example shows the necessity of the condition of the theorem 7.

{\bf Example.} Consider in the space $L^2(T^2)$ (where $T^2=[0,2\pi]^2$)
operator of the form (\ref{1}) with $n=1$, $\nu=2$, $U(x)=\cos x_1$,
$$
\k=\sum_{k\geq 1}c_k\cos kx_2\cos ky_2,\ x=(x_1,x_2), \ y=(y_1,y_2),
$$
where
$$
c^{-1}_k=\int_{T^1}(\cos x_1+1+e^{-k})^{-1}dx_1.
$$
It is clear that these functions $U(x)$ and $\k$ are analytic on $T^2$ and
on $T^2\times T^2$ respectively and the extremal points of the function
$U(x)$ are not isolated. It can be easily checked that the essential
spectrum of the operator $H$ coinsides with the segment $[-1,1]$ and the
numbers $\lambda_n=-1-e^{-n}$, $n=1,2,\ldots$ are eigenvalues of the
operator $H$ lying outside its essential spectrum.

\noindent
{\bf 2.}\ Let $\f\in C^s(\t)$ with natural $s$.
Denote by $J^s_\alpha\f,\ \alpha\in\t$ a $s$-jet of the function $\f$ at
the point $\alpha$ (cf.~\cite{9}). The following theorem is valid.

{\bf Theorem 8.}\ {\it Let $\f(x)$ be a real-valued function from the class
$C^{\mu+3}(\t)$ where $\mu$ is a natural number and let $\alpha\in\t$ be an
extremal point of the function $\f$. If there exists a smooth function
$\psi(x)$ on $\t$ for which the point $\alpha$ is a critical point of
multiplicity $n\leq \mu$ (cf.~\cite{9}) and if $J^{\mu+1}_\alpha\f=
J^{\mu+1}_\alpha\psi$, then $m(\alpha)\leq n+1$.}

{\bf Proof.}\ As the statement of the theorem has a local character so it
is sufficient to prove it in the neighbourhood of zero in ${\bf R}^\nu$.
Without any loss of generality we can suppose that $\alpha=0$ is the
point of minimum of the function $\f$ and $\f(0)=0$. As it is shown
in~\cite{12} the number $n$ is odd.\\
Let $p(x)$ be a $(\mu+1)$-jet of the function $\f$ at the point $\alpha=0$.
It follows from the Tujron theorem~\cite{12} that the function $p(x)$ has a
local minimum at zero.\\
Consider now a one-parameter deformation of the function $p$ of the form
$F_\e(x)=p(x)-\e x^{n+1}_1$ (where $x=(x_1,x_2,\ldots,x_\nu),\ \e>0$.
As the multiplicity of the critical point $\alpha=0$ is equal to $n$, so
$x_1^n\in I_{\nabla p}$ (cf. \cite{9}) where $I_{\nabla p}$ is the local
gradient ideal of the function $p$ at zero. Henceforth there exist such
smooth functions $\{h_k(x)\}$ that in some neighbourhood of zero the
following equality is valid:
$$ x^n_1=h_1(x)\frac{\partial p}{\partial x_1} +
h_2(x)\frac{\partial p}{\partial x_2} + \ldots +
h_\nu(x)\frac{\partial p}{\partial x_\nu}. $$
Hence if $\e$ is small enough we have $I_{\nabla p}=I_{\nabla F_\e}$.
So the map $\nabla F_\e$ at the point $\alpha=0$ has a zero of the
multiplicity $n$ (\cite{9}).

As $\alpha=0$ is the critical point of finite multiplicity for the polynom
$p(x)$, so there exists a ball neighbourhood $V\in{\bf C}^\nu$ of zero such
that for all $x\in\partial V$ (where $\partial V$ is a boundary of $V$) the
inequality $\vert\nabla p(x)\vert\geq\delta>0$.\\
Take $\e$ to satisfy the following inequalities:
\begin{itemize}
\item
$F_\e(x)>0$ for all $x\in\partial V\cap {\bf R}^\nu$,
\item
$\vert\nabla p(x)\vert >\e\vert x^n_1\vert$ for all $x\in\partial V$.
\end{itemize}
By the multidimensional Rouchet theorem (\cite{10}) the maps $\nabla
F_\e(x)$ and $\nabla p(x)$ have equal number of zeros
(counted with multiplicities) in $V$. Therefore both functions $\nabla
F_\e(x)$ and $\nabla p(x)$ have in $V$ a unique zero of multiplicity $n$ at
the point $\alpha=0$, in particular the function $F_\e(x)$ has no real
critical points in $V$. \\
Therefore for any $x\in \bar V\cap {\bf R}^\nu$ and small enough positive
$\e$ we have $F_\e(x)\geq 0$, i.e.
$$ p(x)\geq \e\vert x_1\vert^{n+1}. $$
By the same way we can prove that there exists a neighbourhood $W$ of zero
and a positive number $\e$ such that  the estimates
$$ p(x)\geq\e\vert x_k\vert^{n+1},\quad (k=1,2,\ldots,\nu) $$
hold for all $x\in W$ and it proves the theorem.

{\bf Remark.}\ Validity of the theorem 8 for so-called extremally
non-degenerate polynoms follows from the theorem 1.5 of~\cite{12}.

{}From the theorems 4 and 8 we obtain the following

{\bf Theorem 9.}\ {\it Let $U(x)\in C^{\mu+3}(\t\times\t,\m)$ with some
natural $\mu$. If for any $z\in\Gamma$ and for any $\alpha\in \f^{-1}_z(0)$
where $\f_z(x)=\Delta(x,z)$ there exists a smooth function
$\psi_{z,\alpha}(x)$ on $\t$ for which the point $\alpha$ is critical of
multiplicity $n(\alpha)$ and
$$ J_\alpha^{\mu+1}\f_z=J_\alpha^{\mu+1}\psi_{z,\alpha},\quad \sum_\alpha
n(\alpha)\leq\mu, $$
then for any function $\k\in C^{2\mu+2-\nu/2+0}(\t\times\t,\m)$ the
operator $H$ has
finite number of
eigenvalues
(counted with
multiplicities)
outside the essential spectrum. }

{\bf Acknowledgement.}\ This research was supported by Uzbek Science
Foundation (grant N 44). The second author thanks for partial support by
the International Science Foundation (Soros) (grant N MGM000). We are
grateful to S.~N.~Lakaev and A.~S.~Mishchenko for helpful discussions.

\vspace{2.5cm}
\noindent
Ikramov I.~A.,\quad Sharipov F.\\
Dept. of Mathematics \\
Samarkand State University \\
Samarkand, 703004, Uzbekistan

\end{document}